%
%
 \magnification \magstep1
\font\bold=cmbx10 at 14pt
\font\small=cmr8
\centerline{\bold The homotopy Lie algebra of a complex hyperplane}
\medskip 
\centerline{\bold  arrangement is not necessarily finitely presented.}
\bigskip

\centerline{ Jan-Erik Roos}
\centerline{Department of Mathematics}
\centerline{Stockholm University}
\centerline{SE--106 91 Stockholm, SWEDEN}
\centerline{ e-mail: {\tt jeroos@math.su.se}}
\medskip

\centerline{October 3, 2006}
\medskip
\centerline{Updated: February 9, 2007}
\bigskip

\def\mysec#1{\bigskip\centerline{\bf #1}\nobreak\par}

\def\cite#1{~[{\bf #1}]}
{\openup 1\jot

\mysec{ Abstract.}
\bigskip
We present a theory that produces several examples where the homotopy Lie algebra of a complex
hyperplane arrangement is not finitely presented. We also present examples of hyperplane 
arrangements where the enveloping algebra of this Lie algebra has an irrational Hilbert series.
This answers two questions of Denham and Suciu.


\mysec{0. Introduction}
\bigskip
 Let ${\cal A}=\{H\}$ be a finite set of complex hyperplanes in ${\bf C}^n$,
i.e. a {\it complex hyperplane arrangement} in ${\bf C}^n$
and let $X$ be the complement of their union in ${\bf C}^n$:
$$
X = {\bf C}^n \setminus \bigcup_{H \in {\cal A}}H
$$
The cohomology of $X$ is called the Orlik-Solomon algebra and the
Yoneda Ext-algebra of $H^*(X)$ is a Hopf algebra which is the
enveloping algebra of a graded Lie algebra, which is called the
homotopy Lie algebra of the arrangement ${\cal A}$.
In this paper we calculate explicitly this Lie algebra in several
cases and in particular we show by explicit examples that this Lie algebra
 is not
 necessarily finitely presented and not even finitely generated. Furthermore we present examples
of hyperplane arrangements where the enveloping algebra of this Lie algebra has an 
irrational Hilbert series. This solves two open problems from [De-Su], Question 1.7 p. 321. 
We also have some results about how often these two phenomena occur.

\mysec{1. An explicit example}
\bigskip
It is useful to begin with an explicit example.
Let  ${\cal A} =\{x,y,z,x+y,x+z,y+z\}$
be the well-known complex hyperplane arrangement, which is the smallest formal
arrangement whose Orlik-Solomon algebra is nonquadratic
(cf. [Sh-Yu], p. 487, Example 5.1). We know that the Orlik-Solomon algebra
of ${\cal A}$ is the quotient of the exterior algebra in 6 variables
$e_1,e_2,e_3,e_4,e_5,e_6$  by the twosided ideal generated by the four elements:
$$
(e_2-e_6)(e_3-e_6),\quad (e_1-e_3)(e_3-e_5),\quad (e_1-e_4)(e_2-e_4), \quad
(e_3-e_4)(e_4-e_6)(e_5-e_6)
$$
Thus if we introduce new variables $x_1,x_2,x_3,x_4.x_5,z$ by
$x_i=e_i-e_6$ for $1\leq i \leq 5$ and $z = e_6$ our Orlik-Solomon algebra
can be written (we use that $e_1-e_4 = (e_1-e_6)-(e_4-e_6)=x_1-x_4$ etc.) 
 as a quotient of an exterior algebra:
$$
OS_{\cal A} = {E(x_1,x_2,x_3,x_4,x_5,z) \over (x_2x_3,(x_1-x_3)(x_3-x_5),(x_1-x_4)(x_2-x_4),(x_3-x_4)x_4x_5)}
$$
where $z$ does not occur among the relations.
Therefore the Orlik-Solomon algebra decomposes into a tensor product of algebras (all algebras
are considered over a field $k$ of characteristic 0):
$$
OS_{\cal A} = {E(x_1,x_2,x_3,x_4,x_5)\over (x_2x_3,x_1x_3-x_1x_5+x_3x_5,x_1x_2-x_1x_4+x_2x_4,x_3x_4x_5)}\otimes_k E(z)
$$
where $E(z)$ is the exterior algebra in one variable and where we have used that $x_i^2 = 0$.
Thus the Yoneda Ext-algebra of the Orlik-Solomon algebra is the tensor product of the
Ext-algebra $Ext^*_R(k,k)$ of
$$
R = {E(x_1,x_2,x_3,x_4,x_5)\over (x_2x_3,x_1x_3-x_1x_5+x_3x_5,x_1x_2-x_1x_4+x_2x_4,x_3x_4x_5)}
 \leqno(1.1)
$$ 
and the Ext-algebra $Ext^*_{E(z)}(k,k)= k[Z]$ where the last algebra is the polynomial algebra
in one variable $Z$, dual to $z$.
The last algebra is "innocent" and it therefore follows that the Yoneda Ext-algebra of the
Orlik-Solomon algebra is finitely presented if and only if $Ext^*_R(k,k)$ is so, where $R$ is 
given by (1.1).
But the following automorphism of $R$
$$
x_1 \rightarrow x_1+x_3,x_2\rightarrow -x_2+x_3,x_3 \rightarrow x_3,x_4\rightarrow -x_2+x_3+x_4,x_5\rightarrow x_5+x_3
$$
 transforms  $R$ into the isomorphic algebra
$$
{E(x_1,x_2,x_3,x_4,x_5)\over (x_2x_3,x_1x_5,(x_1+x_2)x_4,x_3x_4x_5)} \leqno(1.2)
$$
which we will still denote by $R$.
But this algebra (1.2) can now be easily analyzed: 
It is the "trivial extension" of a Koszul algebra
$$
S = {E(x_1,x_2,x_3,x_5)\over (x_2x_3,x_1x_5)} \leqno(1.3)
$$
by the following cyclic module $M$ over $S$ :
$$
M = {S \over (x_1+x_2,x_3x_5)}
$$
Recall that the trivial extension of any ring $\Lambda$ by any two-sided  $\Lambda$-module
$N$ is denoted by $\Lambda \propto N$ and consists of the pairs $(\lambda,n)$ with
 $\lambda \in \Lambda$ and $n \in N$ with pairwise addition and multiplication
$$
(\lambda,n).(\lambda',n') = (\lambda .\lambda',\lambda n'+n \lambda')
$$
The Ext-algebra of $R = S \propto M$ can now be analyzed ( cf. e.g. 
[L\"o2,Theorem 3, p. 310-311]):
We have a split extension of Hopf algebras
$$
 k \rightarrow T(s^{-1}Ext^*_S(M,k))\rightarrow Ext^*_R(k,k)\rightarrow Ext^*_S(k,k)\rightarrow k
\leqno(1.4)$$
Here $S$ is the Koszul algebra (1.3) and 
$$
Ext^*_S(k,k)=k< X_1,X_5 > \otimes_k k < X_2,X_3 > 
$$
is the tensor product of two free algebras in the dual variables $X_1,X_5$ and $X_2,X_3$
respectively, and therefore it has global dimension 2. Furthermore
$$
s^{-1}Ext^*_S(M,k)=Ext^{*-1}_S(M,k)
$$ and
 $T(s^{-1}Ext^*_S(M,k))$ is the free algebra on the
graded vector space (for the *-grading in Ext) $s^{-1}Ext^*_S(M,k)$ and has global dimension 1.
The spectral sequence of extensions of Hopf algebras (1.4) [Ro2]
$$
E^2_{p ,q} = Tor_p^{Ext_S^*(k,k)}(k,Tor_q^{T(s^{-1}Ext^*_S(M,k))}(k,k) \Rightarrow
Tor_n^{Ext^*_R(k,k)}(k,k) \,( \,= H_n) \leqno (1.5)
$$
shows immediately that $Ext^*_R(k,k)$ har global dimension 3.
Furthermore (1.5) degenerates into a long exact sequence:
$$
0 \rightarrow E^2_{2,1}\rightarrow H_3 \rightarrow E^2_{3,0}\rightarrow 
E^2_{1,1}\rightarrow   H_2 \rightarrow E^2_{2,0}\rightarrow E^2_{0,1}\rightarrow 
H_1\rightarrow E^2_{1,0} \rightarrow 0 \leqno(1.6)
$$
where the natural maps $H_i \longrightarrow E^2_{i,0}$ are {\it onto}. Indeed,
the natural ring projection map $S \propto M \longrightarrow S$
is split by the natural ring inclusion $S \longrightarrow S \propto M$ and
this leads to a splitting on the Ext-algebra level.
Thus we have exact sequences
$$
0 \longrightarrow E^2_{i-1,1} \longrightarrow H_i \longrightarrow E^2_{i,0 }
\longrightarrow 0.
$$
 Now for any graded connected algebra $A$ over k, $Tor_1^A(k,k)$ measures the minimal number
of generators of $A$, $Tor_2^A(k,k)$ measures the minimal number of relations between these
generators, $Tor_3^A(k,k)$ the minimal number of relations between these relations etc.,
cf. the Chapter 1 ("Pr\'esentations d'alg\`ebres connexes") of Lemaire [Lem1].
Therefore $H_1$ in (1.6) measures the minimal number of generators
of the Ext-algebra $Ext^*_R(k,k)$ and therefore $H_1$ is finite-dimensional if and
only if the Ext-algebra $Ext^*_R(k,k)$ is finitely generated.
Similarly $H_2$ measures the minimal number of relations in a minimal presentation of
$Ext^*_R(k,k)$ and $H_3$ measures the minimal number of relations between these
relations. Since the $E^2_{i,0}$ are all finite-dimensional we
are led to the study of
$$
E^2_{i,1} =Tor_i^{Ext^*_S(k,k)}(k, s^{-1}Ext^*_S(M,k)), \, {\rm for} \,i\geq 1\leqno(1.7)
$$
where the left $Ext^*_S(k,k)$-module structure of $s^{-1}Ext^*_S(M,k)$ is given
 by the Yoneda product
({\it cf.} again Theorem 3, p. 310-311 of L\"ofwall [L\"o2]). Note that underlying our spectral
 sequence is the Hochschild-Serre
spectral sequence and that we are in the skew-commutative setting, whereas 
[L\"o2] is in the commutative case, but similar
(easier) proofs work here in our case.
Thus to show that $Ext^*_R(k,k)$ is not finitely generated we have to show that $H_1$ is infinite-dimensional,
i.e. that
 $E^2_{0,1}$ is so, i.e. (cf. (1.7)) that $s^{-1}Ext^*_S(M,k)$ 
needs an infinite number of generators as an $Ext^*_S(k,k)$-module, i.e. we have to 
study the $S$-resolutions of $M = S/(x_1+x_2,x_3x_5)$. We also need the extra grading on $R$
and $S$ so that we should indeed write $R = S \propto s^{-1}M$.
Now we denote the $S$-ideal $(x_1+x_2,x_3x_5)$ by $I$ so that $M = S/I$.
First we observe that if we apply the functor $Ext_S^*(.,k)$ to the exact sequence of 
graded left $S$-modules
$$
0 \longrightarrow I \longrightarrow S \longrightarrow S/I \longrightarrow 0 \leqno(1.8)
$$
we obtain the isomorphisms of left $Ext_S^*(k,k)$-modules:
$$
Ext^{*-1,t}_S(I,k)\buildrel{ \sim} \over \longrightarrow Ext^{*,t}_S(S/I,k),\,\, 
{\rm for}\,\, * \geq 1 \leqno (1.9) 
$$
where we also have inserted the inner grading $t$ that comes from the fact that 
(1.8) is an exact sequence of graded modules. Note that $S$ is a Koszul algebra, so that
only the $Ext^{i,i}_S(k,k)$ are different from zero and we still denote them by $Ext^{i}_S(k,k)$.
Next we note that the two ideals
  $I_1=(x_1+x_2)$ and $I_2 =(x_3x_5)$ in $S$ have zero intersection.
Therefore $I = I_1 \oplus I_2$ and the $Ext_S^*(k,k)$-module to the left in (1.9) decomposes into
a direct sum of $Ext_S^*(k,k)$-modules:
$$
Ext^{*-1,t}_S((x_1+x_2),k) \oplus Ext^{*-1,t}_S((x_3x_5),k) \leqno(1.10)
$$
But $x_3x_5$ is in the socle of $S$ and therefore we have as graded $S$-modules that
$(x_3x_5) \buildrel{ \sim} \over \longrightarrow s^{-2}k$ so that the right summand of (1.10)
is isomorphic to $Ext^{*-1,t}_S(s^{-2}k,k)$, i.e. to $Ext^{*-1,t-2}_S(k,k)$.
It remains to analyze the left summand of (1.10).
But it is easy to see that $Ann_S((x_1+x_2))=I$ so that
the graded sequence of $S$-modules
$$
0 \longrightarrow s^{-1}I \longrightarrow s^{-1}S \buildrel{ (x_1+x_2)} \over \longrightarrow S \leqno (1.11)
$$
where we multiply to the right with $x_1+x_2$ is exact. Therefore we have a short exact sequence
of graded left $S$-modules:
$$
0 \longrightarrow s^{-1}I \longrightarrow s^{-1}S \longrightarrow (x_1+x_2)
 \longrightarrow 0 \leqno(1.12)
$$
leading to the isomorphism of left $Ext_S^*(k,k)$-modules:
$$
Ext^{*-1,t}_S(s^{-1}I,k) \buildrel{ \sim} \over \longrightarrow Ext^{*,t}_S((x_1+x_2),k)
 ,\,\, {\rm for}\,\, * \geq 1 \leqno (1.13)
$$
and using (1.9) once more we obtain that
$$
Ext^{*-1,t}_S(s^{-1}I,k) = Ext^{*-1,t-1}_S(I,k)\buildrel{ \sim} \over \longrightarrow
Ext^{*,t-1}_S(S/I,k) \leqno (1.14)
$$
leading to the final isomorphism (combining (1.9),(1.10),(1.13),(1.14)):
$$
Ext^{*,t}_S(S/I,k)\buildrel{ \sim}\over \longrightarrow Ext^{*-1,t-1}_S(S/I,k) \oplus
Ext^{*-1,t-2}_S(k,k) \leqno(1.15)
$$
for $* \geq 1$, where the summand $Ext^{*-1,t-2}_S(k,k)$ is non-zero only if $*=t-1$.
This proves everything since we see, using (1.15) that
$Ext^{*,t}_S(S/I,k)$ needs a new $Ext^*_S(k,k)$-generator for $*=t-1$ for each $t=2,3,4,\ldots$.

In particular, if we introduce for any graded module $N$ over a graded $k$-algebra $G$
 the double series
$$
P_G^N(x,y) = \sum_{i\geq 0,j\geq 0} |{Ext}_G^{i,j}(N,k)|x^iy^j  \leqno(1.16)
$$
(where as always for a $k$-vector space $V$ we denote by $|V|$
its dimension) and if we denote $P_G^k(x,y)$ by $P_G(x,y)$,
we then deduce from (1.15) and the fact that $P_S(x,y)=1/(1-2xy)^2$ that
$$
P_S^{S/I}(x,y) = {1 \over 1-xy}+ {xy^2 \over (1-xy)(1-2xy)^2}
$$
so that 
$$
P_{S\propto s^{-1}S/I}(x,y) = {P_S(x,y) \over 1-xyP_S^{S/I}(x,y)}
$$
leading to the theorem:

{\bf THEOREM 1.1.} The  Orlik-Solomon algebra  of the
complex hyperplane arrangement ${\cal A}= \{x,y,z,x+y,x+z,y+z\}$ is the tensor product
of the exterior algebra
in one variable with an algebra $R$ whose Yoneda Ext-algebra $Ext^*_R(k,k)$
has a bigraded generating series:
$$
P_R(x,y)={P_S(x,y) \over 1-xyP_S^M(x,y)} = {1-xy \over 1-6xy+12x^2y^2-x^2y^3-8x^3y^3} \leqno(1.17)
$$
where $S$ and $M$ are defined above.
Furthermore, the Ext-algebra $Ext^*_R(k,k)$ has global dimension 3 and it 
has 5 generators in degree 1 and needs one new generator in each degree $\geq 2$.
 In particular the homotopy Lie algebra of ${\cal A}$
is not finitely generated.

We will come back to this result in the next section.
\bigskip
\mysec{2. The holonomy and homotopy Lie algebra of an arrangement.}
\medskip
The preceding analysis of the ${\cal A}$ arrangement in section 1
was intended to give the "simplest possible proof" that the
Ext-algebra $Ext^*_R(k,k)$ was not finitely generated.
 However, in order to be able to analyze more cases we need
a more general theory. I will here briefly describe the basics of such a theory and apply 
it as an alternative to our
first case and then treat another case of arrangements where we can prove that the 
homotopy Lie algebra is also non-finitely presented.
Note that the graded algebra $R$ of the previous section has Hilbert Series $1+5z+7z^2$. 
Let us now start with any algebra $R$ which is a quotient of an exterior algebra 
$E(x_1,x_2,\ldots,x_n)$ by a homogeneous ideal $J$ generated by
elements of deegree $\geq 2$. Thus $R = E(x_1,\ldots,x_n)/J$.  Let $m$ be the ideal of
 $R$ generated by $(x_1,\ldots,x_n)$,
and consider the exact sequence of left $R$-modules:
$$ 
0 \longrightarrow m/m^2 \longrightarrow R/m^2 \longrightarrow R/m \longrightarrow 0, \leqno(2.1)
$$
Now apply the functor $Ext^*_R( ,k)$ to the exact sequence (2.1). We get a long exact sequence
which can be written as an 
exact sequence of left $Ext^*_R(k,k)$-modules (we use the Yoneda product):
$$
0 \rightarrow s^{-1}{\overline S_m} \rightarrow s^{-1}{\overline {Ext}}^*_R(R/m^2,k)
 \rightarrow Ext^*_R(k,k)\otimes Ext^1_R(k,k)
\rightarrow Ext^*_R(k,k) \rightarrow S_m  \rightarrow 0 \leqno(2.2)
$$
where $S_m$ is defined as the image of the natural map
$$
Ext^*_R(k,k) \longrightarrow Ext^*_R(R/m^2,k)  \leqno(2.3)
$$
and where e.g. $\overline S_m$ means that we take the elements of $S_m$ of  degrees $> 0$
and where $s^{-1}$ is "the suspension" as before.
Now the Ext-algebra $B = Ext^*_R(k,k)$ is bigraded and we recall that its bigraded Hilbert
 series is denoted
by
 $$
P_R(x,y) = \sum_{i,j \geq 0} |Ext^{i,j}_R(k,k)|x^iy^j = B(x,y),
$$
where as always for a $k$-vector space we denote by $|V|$ its dimension.
The subalgebra $A$ of $Ext^*_R(k,k)$ generated by $Ext^1_R(k,k)$ is also bigraded but it is situated
on the diagonal so that the corresponding bigraded Hilbert series 
$$
A(x,y) = A(xy,1)  \buildrel{ def} \over = A(xy)
$$
Now take the alternating sum of the two-variable Hilbert series of (2.2). We obtain:
$$
S_m(x,y)-B(x,y)+xyB(x,y)|m/m^2|-x{\overline P}_R^{R/m^2}(x,y)+x(S_m(x,y)-1) = 0 \leqno(2.4)
$$
where
 $$
{\overline P}_R^{R/m^2}(x,y) = \sum_{i>0,j\geq i} |Ext_R^{i,j}(R/m^2,k)|x^iy^j \leqno (2.5)
$$
We now make three fundamental observations:
\medskip
1) $A$ is a sub Hopf algebra of $B$ and therefore according to a result of Milnor and Moore [Mi-Mo]
$B$ is free over $A$. Thus $S_m = B \otimes_A k$ has bigraded Hilbert series
$$
S_m(x,y) = B(x,y)/A(xy) \leqno (2.6)
$$
\medskip
2) If $m^3 = 0$ we have an isomorphism
of left $Ext^*_R(k,k)$-modules
$$
{\overline {Ext}}^*_R(R/m^2,k) \simeq Ext^*_R(k,k)\otimes Ext^1_R(R/m^2,k) \leqno(2.7)
$$
so that 
$$
{\overline P}_R^{R/m^2}(x,y) = B(x,y)xy^2 |m^2/m^3| \leqno(2.8)
$$
Therefore the equality (2.4) can be written
$$ 
{B(x,y)\over A(xy)}= B(x,y)-xyB(x,y)|m/m^2|+x^2y^2B(x,y)|m^2/m^3|-x({B(x,y)\over A(xy)}-1) \leqno(2.9)
$$
which is another way of writing (divide by $xB(x,y)$ and use the 
notation $R(z) = 1+|m/m^2|z+|m^2/m^3|z^2$ for the Hilbert series of $R$ )
$$
1/B(x,y) = (1+1/x)/A(xy)-R(-xy)/x \leqno(2.10)
$$
which is a formula due to L\"ofwall [L\"o1].

\medskip
3) In the case when $m^3=0$ the 3 middle terms of (2.2) are free $Ext^*_R(k,k)$-modules
so that $s^{-1}{\overline S}_m$ is a third syzygy of a minimal graded $Ext^*_R(k,k)$-resolution of
$S_m$. We therefore obtain the isomorphism:
$$
Tor_{i,*}^B(k,S_m)  \simeq  Tor_{i-3,*}^B(k,s^{-1}{\overline S}_m) 
=  Tor_{i-3,*-1}^B(k,{\overline S}_m) 
\quad {\rm for}\, i\geq 3 \leqno(2.11)
$$
Now apply $Tor_{i}^B(k,)$ to the exact sequence:
$$
0 \longrightarrow {\overline S}_m \longrightarrow S_m \longrightarrow k \longrightarrow 0 \leqno(2.12)
$$
We get a long exact sequence:
$$
\ldots\rightarrow Tor_{n+1}^B(k,k) \rightarrow Tor_n^B(k,{\overline S}_m)
 \rightarrow Tor_n^B(k,S_m) \buildrel{ \varphi_n} \over \rightarrow Tor_n^B(k,k)\rightarrow Tor_{n-1}^B(k,{\overline S}_m)\ldots\leqno(2.13)
$$
Furthermore, since $B$ is $A$-flat
$$
Tor_n^B(k,S_m) = Tor_n^B(k,B\otimes_A k) = Tor_n^A(k,k)
$$
and $\varphi_n : Tor_n^A(k,k) \rightarrow Tor_n^B(k,k)$ is induced by the natural inclusion
$A \rightarrow B$ which is split by a ring map in the other direction:
Divide  $B = Ext_R^{*,*}(k,k)$ by the twosided ideal generated by $\oplus_{j>i>0}Ext^{i,j}_R(k,k)$.
Thus the maps $\varphi_n$ in (R) are monomorphisms and (2.13) splits into short exact sequences,
using (2.11):
$$
 0 \longrightarrow Tor_{i,j}^A(k,k)\longrightarrow Tor_{i,j}^B(k,k) \longrightarrow
Tor_{i+2,j+1}^A(k,k) \longrightarrow 0 \leqno(2.14)
$$ 

Now we can summarize:

{\bf THEOREM 2.1.} Let $R$ be a quotient of an exterior algebra (finite number of generators 
in degree 1) by a homogeneous ideal generated by elements of degree $\geq 2$. Let $m$ be the augmentation ideal 
of $R$. Assume that $m^3=0$ . Let $B = Ext^*_R(k,k)$ be the Yoneda Ext-algebra and let
$A$ be the subalgebra of $B$, generated by $Ext^1_R(k,k)$. Then the exact sequences
(2.14) hold.
In particular

a) $B$ is finitely generated if and only if the graded vector space $Tor_{3,*}^A(k,k)$
 has finite dimension.

b) $B$ is finitely presented  if and only if the graded vector spaces $Tor_{3,*}^A(k,k)$
and $Tor_{4,*}^A(k,k)$
 have finite dimension.

c) $B$ is finitely presented and has a finite number of relations between the minimal relations
if and only if 
the graded vector spaces $Tor_{3,*}^A(k,k)$ , $Tor_{4,*}^A(k,k)$ and $Tor_{5,*}^A(k,k)$
have finite dimension.

Etc.

Note that $A$ is the enveloping algebra of a graded Lie algebra (the holonomy
Lie algebra) whose ranks are equal to the ranks of the lower central series (LCS) of the
fundamental group of the hyperplane complement (cf. section 6 below).
Note also that in general the Hilbert series $A(x)$ of $A$ when $m^3=0$ is obtained from (2.10):
replace $x$ by $x/y$ in that formula and  put $y=0$. This gives
$$
1/P_{R}(x/y,y)|_{y=0} = 1/A(x)
$$
We can get an alternative proof of assertion about generators of the $Ext$-algebra 
in Theorem 1.1 above, using only the formula (1.17) there. The preceding recipe gives in that case:
that $A(x) = (1-x)/(1-2x)^3$. Now recall that for any graded algebra $A$ we have the 
following formula for the relation
between its Hilbert
series $A(z)$ and the Hilbert series $Tor_{i,*}^A(k,k)(z)$ of the graded Tor:
(cf. eg. Lemaire [Lem1, Appendix A2]):
$$
{1\over A(x)} = \sum_{i \geq 0}(-1)^iTor_{i,*}^A(k,k)(x) \leqno(2.15)
$$
Since $A$ in the case of Theorem 1.1 has global dimension 3 and 5 generators in degree 1 and
7 relations in degree 2, (2.15) gives that $Tor_{3,*}^A(k,k)(x) = x^3/(1-x)$
 so that we see once more that
$Tor_{3,i}^A(k,k)$ is one-dimensional for all $i\geq 3$.

\noindent {\it Remark 2.2.} If $m^3 \neq 0$ but more generally
$$
{\overline {Ext}}^*_R(R/m^i,k) \longrightarrow {\overline {Ext}}^*_R(R/m^{i+1},k) \leqno(2.16)
$$
is 0 for $i\geq 2$ then we have the same conclusion as in Theorem 2.1 but the proof is slightly
different, since now $N={\overline {Ext}}^*_R(R/m^2,k)$ is not free as a $B=Ext^*_R(k,k)$-module but it has a finite
homological dimension and the corresponding $Tor_i^B(k,N) \simeq m^{i+2}/m^{i+3}$ are finite-dimensional.
The condition (2.16) is sometimes, but not always satisfied  if $m^4 =0$, but in the last case
one can prove the validity of the fomula (2.10) is {\it equivalent} to the assertion
that the map (2.16) is zero for $i\geq 2$ (only the case $i=2$ is of course important in this case).
This will be used below when we study graphic arrangements.
Furthermore the important formula (2.10) which we will now write as 
$$
1/P_R(x,y) = (1+1/x)/R^!(xy) -R(-xy)/x \leqno(2.17)
$$
 is still valid under (2.16) but here the Hilbert series
$R(z)$ might be a polynomial of degree $> 2$.
Note that we have written $R^!$ (instead of $A$); it is
the Koszul dual of $R$.
\bigskip
\noindent {\it Remark 2.3.} 
The formula (2.17) in Remark 2.2 is the special case
 is a special case ($n = 3$) of a whole family of formulae ($n \geq 3$)
$$
{1\over P_R(x,y)} = {(1-(-x)^{2-n})\over R^!(xy)}+ R(-xy)(-x)^{2-n} \leqno(2.18_n)
$$
The validity of $(2.18_n)$ is a consequence of the fact 
that the so-called Koszul complex
$R^! \otimes_k Hom_k(R,k)$ has only non-zero homology groups in degree 0 and in degree $n-1$.
For more details about this, cf. the appendix B (Theorem B.4) by L\"ofwall [L\"o4] 
to [Ro3] and also [Ro5].
We will say here that $R$ satisfies  $L_n$ if $(2.18_n$) holds true.
For Orlik-Solomon algebras with $m^4=0$ we still have the formula (2.17) since the algebra
is the tensor product of
an algebra with $m^3=0$ and an ``innocent'' algebra $E[z]$.
In the section 5 where we study the case of graphic arrangements we will see that for any $n \geq 3$
there are examples where the condition $L_n$ holds (namely the Orlik-Solomon algebra of 
the graphic arrangement corresponding
to an $n+1$-gon for $n\geq 3$), but that there are also examples where none of these
conditions is satisfied (these cases can however sometimes be handled with the method of [Ro5]).

\mysec{3. Some other hyperplane arrangements}
\bigskip
Some of the cases from [Su] can be treated in the same way as in 
section 2. Here we just briefly describe the results for the so-called $X_2$-arrangement
 which is defined by the 
polynomial $xyz(x+y)(x-z)(y-z)(x+y-2z)$.
Now the Orlik-Solomon algebra can be written
as a quotient of the exterior algebra in 7 variables $E(e_1,e_2,e_3,e_4,e_5,e_6,e_7)$
by the ideal generated by  five elements:
$$
(e_5-e_7)(e_6-e_7),(e_3-e_7)(e_4-e_7),(e_2-e_6)(e_3-e_6),(e_1-e_5)(e_3-e_5),(e_1-e_2)(e_2-e_4)
$$
Now isolate $e_3$, i.e. introduce variables $a=e_1-e_3,b=e_2-e_3,c=e_4-e_3,d=e_5-e_3,e=e_6-e_3,f=e_7-e_3$.
Now the relations in the Orlik-Solomon algebra do not contain $e_3$ and this algebra is now a tensor algebra
of the quotient:
$$
R={E(a,b,c,d,e,f)\over (ab-ac+bc,ad,be,cf,de-df+ef)}
$$
with the exterior algebra in one variable $z=e_3$. Therefore we are lead to the analysis 
 of the Yoneda Ext-algebra of
the quotient $R$ above 
whose Hilbert series is $R(t)=1+6t+10t^2$.
Furthermore, let
$$
S={E(a,b,c,d,e,f)\over (ab-ac+bc,ad,be,cf)}
$$
and consider the ring map:
$$
S \longrightarrow S/(de-df+ef) = R  \leqno(3.1)
$$
It is not difficult to show that (3.1) is a so-called Golod map (cf. [Le2] and the
litterature cited there). One finds that
$S^!(t) = (1-t)^2/(1-2t)^4$ ,that 
 $R^!(t) = (1-t)^4/(1-2t)^5$ and more precisely
that
$R^!$ has global dimension 5 and that:
$$
\sum_{i \geq 0}|Tor_{3,i}^{R^!}(k,k)|z^i = 5z^4+{6z^5\over(1-z)}\, , \leqno(3.2)
$$
$$
\sum_{i \geq 0}|Tor_{4,i}^{R^!}(k,k)|z^i = 2z^6+{(6-z)z^7\over(1-z)^2}  \leqno(3.3)
$$
and
$$
\sum_{i \geq 0}|Tor_{5,i}^{R^!}(k,k)|z^i = {z^{10}\over(1-z)^4} \leqno(3.4)
$$
Thus using the theory from section 2 we see that the homotopy Lie algebra of the
arrangement $X_2$ is "extremely non-finitely presented":
It needs an infinite number of generators (3.2), 
and the minimal number of relations between a minimal
 system of generators is infinite (3.3) and the minimal number of relations between the relations is
infinite (3.4).
Among the graphic arrangements (cf. section 5 below for more details) there are however more
arrangements with finitely presented $Ext$-algebra than infinitely presented ones.

We finish this section with one unsolved case:
Recall that the so-called non-Fano arrangement is the hyperplane arrangement
defined by $xyz(x-y)(x-z)(y-z)(x+y-z)$.
In this case the corresponding $R$ (we have eliminated one variable as above) has Hilbert series
$(1+3t)^2$ but the corresponding $R^!(t)$ is rather complicated.
We have however managed to calculate the LCS-ranks two steps higher than in Sucio [Su],
using the Backelin et al programme BERGMAN [B]; with the notations of [Su] we have
$\phi_8 = 3148$ and $\phi_9 = 9857$, but for the last result we needed 64-bits $PSL$ on an AMD
opteron machine with 12 GB of internal memory.

\mysec {4. Arrangements with irrational Hilbert series.}
In [De-Su] it is also asked if the enveloping algebra of the homotopy Lie algebra of an arrangement
 can have an irrational Hilbert series.

We will here describe one case we have found where this is conjecturally true 
and a second case where this is {\it proved} to be true. This development is rather recent:
we found the second case only recently and the first (more complicated) case is the well-known
Mac Lane arrangement, whose amazing homological properties we also discovered recently.
The proof in the second case (the first case is probably treated in a similar but more complicated
way) is based on ideas of the present paper, but involves a lot more new ideas and will be 
presented in another paper in preparation [Ro6]. Let us just indicate some more details:

{\bf First case} (the Mac Lane arrangement):

Recall that the Mac Lane arrangement defined by the annihilation of the polynomial
 $Q = xyz(y-x)(z-x)(z+\omega y)(z+\omega^2 x+\omega y)(z-x-\omega^2 y)$  in ${\bf C}^3$
where $\omega=e^{{2\pi i}/3}$.
It is not difficult to see that with the notations of our section 1 above 
the Orlik-Solomon algebra of the Mac Lane arrangement is $R \otimes E[z]$ where
$R$ is the quotient of the exterior algebra $E(x_1,x_2,x_3,x_4,x_5,x_6,x_7)$ in 7 variables
with the ideal generated by the 8 quadratic elements
$x_1x_2-x_1x_4+x_2x_4,x_1x_3-x_1x_5+x_3x_5,x_1x_6-x_1x_7+x_6x_7,\break
x_2x_3-x_2x_6+x_3x_6,x_4x_5-x_4x_7+x_5x_7,x_2x_5,x_4x_6,x_3x_7
$.
This ring $R$ has Hilbert series $R(z)=1+7\,z+13\,z^2$
and therefore the formula (2.17) above can be applied, and the only thing needed
to be proved is that
the Koszul dual $R^!$ of $R$ has an irrational Hilbert series.
But this Koszul dual $R^!=U(g)$ is the quotient of the free associative algebra $k<X_1,X_2,\ldots,X_7>$
in the seven dual variables by the two-sided ideal generated by the 13 dual relations among the
Lie-commutators $[X_i,X_j] = X_iX_j-X_jX_i$ for $i\neq j$:
$[X_1,X_2]+[X_1,X_4],[X_1,X_4]+[X_2,X_4],
[X_1,X_3]+[X_1,X_5],[X_1,X_5]+[X_3,X_5],
[X_1,X_6]+[X_1,X_7],[X_1,X_7]+[X_6,X_7],
[X_2,X_3]+[X_2,X_6],[X_2,X_6]+[X_3,X_6],
[X_4,X_5]+[X_4,X_7],\break [X_4,X_7]+[X_5,X_7],
[X_2,X_7],[X_3,X_4],[X_5,X_6]$
so that the Lie algebra $g$ (it is called the holonomy Lie algebra) in $R^!=U(g)$ is the quotient
of the free Lie algebra in 7 variables
by the ideal generated by the 13 Lie commutators above. Now we have the formula
$$ 
{1\over R^!(z)} = {1\over U(g)(z)} =  \prod_{n=1}^\infty(1-t^n)^{\phi_i}
$$
where the ${\phi_i}$ are the so-called LCS-ranks (the lower central series ranks).
But these ranks can be calculated by a programme by Clas L\"ofwall [L\"o3] which is called 
{\tt liedim.m} and runs under Mathematica. It gives (in characteristic 0) the ranks
$7,8,21,42,87,105,172,264,476,\break 816,1516,2704,5068,9312,17484, \ldots$ but for the
 higher ranks
you need the C-version of the programme unless your computer has lots of internal memory.
Thus the Hilbert series $R^!$ can be calculated in degrees $\leq 15$ and in these degrees 
it is described by the degree $\leq 15$ part of following rather amazing formula:
$$
{1\over R^!(t)}={(1-2t)^8\over (1-t)^9}(1-t^3)^5(1-t^4)^{18}(1-t^5)^{39}(1-t^6)^{33}
\prod_{n=4}^\infty(1-t^{2n-1})^{28}(1-t^{2n})^{24}
$$
We indicate a possible proof below.

{\bf Second case} (a "simplification-degeneration" of the Mac Lane arrangement):
In the Mac Lane arrangement above $\omega$ being a primitive third root of unity
satisfies ${\omega}^2= -\omega -1$ so if you replace ${\omega}^2$ by $-\omega -1$ in
the definition of the Mac Lane arrangement above {\it and then } (this is brutal!)
 put $\omega = 1$ you obtain
a new arrangement $mleas$ in ${\bf C}^3$ defined by the polynomial with integer coefficients:
$$
Q_{eas}=xyz(y-x)(z-x)(z+y)(z-2x+y)(z-x+2y)  \leqno(4.1)
$$
The amazing thing now is that the corresponding hyperplane arrangement has an almost identical
Orlik-Solomon algebra $R_{eas} \otimes E[z] $ but its homological properties are dramatically 
different:
We have indeed that $R_{eas}$ has all the relations of $R$ {\it with the exeption of
the relation  $x_3x_7$ which is replaced by $x_3x_6x_7$}, so that $R = R_{eas}/(x_3x_7)$.
 and the Hilbert series 
$R_{eas}(t) =1+7\,t+14\,t^2$, which is close to $R(t) =1+7\,t+13\,t^2$. But now we have:

{\bf Theorem 4.1.} The Koszul dual $R_{eas}^!$ of the hyperplane arrangement $mleas$ has the following
Hilbert series
$$
{1 \over R_{eas}^!(t)} = {(1-2\,t)^7 \over (1-t)^7}\prod_{n=3}^\infty(1-t^n) \leqno(4.2)
$$
{\bf Corollary 4.2.} The Ext-algebra of the algebra $R_{eas}$ corresponding the hyperplane arrangement
$mleas$ has a transcendental Hilbert series.

Proof of Corollary 4.2: The condition $L_3$ is satisfied since $m^3=0$.

SKETCH OF PROOF OF THEOREM 4.1 --- We have $R_{eas}^! = R^!/([X_3,X_7])$. Furthermore
$R_{eas}^!$ is the enveloping algebra of a Lie algebra $g_{eas}$. Now divide out $g_{eas}$ with
the  Lie element of degree 3 $[X_6,[X_7,X_5]]$. We get an exact sequence of Lie algebras:
$$
0\longrightarrow ker \longrightarrow g_{eas} \longrightarrow g_{eas}/([X_6,[X_7,X_5]]) 
\longrightarrow 0 \leqno(4.3)
$$
where $ker$ is defined by (4.3). I now claim that the Hilbert series of the enveloping algebra of $quot = g_{eas}/([X_6,[X_7,X_5]])$
is $(1-z)^7/(1-2\,z)^7$.  Indeed
the underlying Lie algebra has a basis of 7 elements $X_1,X_2,\ldots,X_7$
in degree 1.  From now on we denote to simplify the elements 
$[X_i,[X_j,[X_k,\ldots ]]]$ by $ijk\ldots$. With these notations we
have the following basis of seven elements in degree 2:
$$
42, \,52, \, 53 ,\,63,\, 64, \,75,\,76 
$$
These seven elements commute pairwise, and we have 14 elements of degree 3.
One can prove that the Lie algebra decomposes in the sense of [Pa-Su], i.e. that the Lie algebra $quot$ in degrees $\geq 2$
is a direct sum of 7  degree $\geq 2$ parts of free Lie algebras in 2 variables.
One can not use [Pa-Su] directly, but if $S$ is $R_{eas}$ without the relation $x_4x_5-x_4x_7+x_5x_7$ then S comes from the so-called
$X_2$ arrangement which decomposes [Pa-Su] and has $S^!(z)=(1-z)^3/(1-2z)^5$ and the map $S \rightarrow R_{eas}$
can be analyzed as in section 3 above.

We now continue analyzing the kernel $ker$ in (4.3).
Clearly $675=[X_6,[X_7,X_5]]$ is in this kernel and so are also the degree 4 elements
$i675=[X_i,[X_6,[X_7,X_5]]]$ for $i=1\ldots 7$. But we only get possibly non-zero elements 
(with different signs) for $i=5$ and $i=6$,
( the last one can be written 6775) and in next degree we similarly get only one element 67775 etc.
 Therefore $ker$ is $\leq 1$-dimensional
in each degree $\geq 3$. Next we have the formulae in $g_{eas}$ 
$$
-[X_5,42] = [X_6,42]=-[X_4,52] = [X_6,52]=[X_2,53] = -[X_6,53] = [X_1,63] = 
$$
$$
 -[X_5,63] =-[X_2,64] = [X_5,64]=-[X_1,75] = [X_6,75]= -[X_4,76] = [X_5,76]=675    \leqno(4.4)
$$
The other commutators lie in $quot$. Therefore if we define a graded vector space of dimension
 1 in each degree $\geq 3$ by
$V_* = ke_3 \oplus ke_4 \ldots $
where the $X_i$ operate as zero, with the exception of $X_5.e_n=e_{n+1}$ and $X_6.e_n=-e_{n+1}$,
 we get a $quot$-module $V$ and by calculating $H^2(quot,V)$ we find that there is a 2-cocycle
$\gamma: quot \times quot \longrightarrow V$ on $quot$ with values in $V$ which in lower degrees starts
as in (4.4). Thus we have a Lie algebra $g_{\gamma}$ defined by this cocycle which sits in the
middle of an extension of Lie algebras where the kernel $V$ is abelian:
$$ 
0\longrightarrow V \longrightarrow g_{\gamma} \longrightarrow quot \longrightarrow 0 \leqno(4.5)
$$
One now uses the Hochschild-Serre spectral sequence of the extension (4.5) to show that
$Tor_{2,*}^{U(g_{\gamma})}(k,k)$ is concentrated in degrees 2 .
Now $g_{eas}$ and $g_{\gamma}$ are isomorphic (they have "the same" generators and relations). 
Thus $ker$ is one-dimensional in each degree $\geq 3$ This gives the formula (4.2)
and Theorem 4.1 follows.

{\it Remark 4.3.} The previous reasoning with explicit cocycles is analogous to, but a little more
complicated than L\"ofwall's and mine version of the Anick solution of the Serre-Kaplansky
irrationality problem (cf. [L\"o-Ro2], [Ro7,pages 454-456], [An]) as well as Lemaire's 
Bourbaki talk about these questions [Lem2]).

Now what about the Mac Lane arrangement ?
The results are similar to those of the easier case just described.
In fact the Lie algebra $g_{eas}$ just studied has no center and the kernel Lie algebra
$ker$ of (4.3) is abelian.
In order to do the similar reasoning for the Mac Lane (ML) arrangement one needs to divide the Lie algebra
$g_{ML}$ with the following {\it five} cubic Lie algebra elements:
$$
[X_5,[X_7,X_3]],[X_6,[X_7,X_3]],[X_7,[X_6,X_3]],[X_7,[X_6,X_4]],[X_6,[X_7,X_5]] \leqno(4.6)
$$
leading again to a quotient Lie algebra $quot_{ML}$ whose enveloping algebra has Hilbert series
$(1-z)^9/(1-2\,z)^8$
We still get an exact sequence of Lie algebras:
$$
0 \longrightarrow ker_{ML} \longrightarrow g_{ML} \longrightarrow quot_{ML}
 \longrightarrow 0 \leqno(4.7)
$$
But now $ker_{ML}$ is generated by the five elements (4.6) and is still situated in degrees $\geq 3$
 where its dimensions are
$$
5,18,39,33,28,24,28,24,28,24,28,24,28,\ldots \leqno(4.8) 
$$ 
and there is still a two-cocycle describing the extension (4.7). However now $g_{ML}$
contains central elements, all situated in $ker_{ML}$ and $[ker_{ML},ker_{ML}]$
is contained in the center, which is 1-dimensional in degree 5, 9-dimensional in degree 6,
and 4-dimensional in degrees $2n+1$ for $n \geq 3$ and zerodimensional in even degrees $\geq 8$.
If we divide out by the center the Lie algebra $ker_{ML}/center$ is still situated in degrees
$\geq 3$ but its dimensions there are now
$$
5,18,38,24,24,24,24,24,24,24,24,24,24,\ldots \leqno(4.9) 
$$ 
Furthermore $ker_{ML}/center$ is abelian and $quot_{ML}$ operates on it in a 
similar but more complicated way than for $g_{eas}$ above. But so far all this has only been proved
 in degrees $\leq 15$.

Although the $g_{eas}$ irrational case is much easier than the $g_{ML}$ case it still uses a 
hyperplane arrangement with 8 hyperplanes
leading to an algebra $R_{eas}^!$ in {\it seven variables}. One might still wonder it it would be
 possible to simplify further, i.e. get a hyperplane arrangement 
with seven or six hyperplanes and having an irrational series:

In [Ro4] we have in particular described all homological possibilities for the quotient of an 
exterior algebra in $\leq 5$ variables
by an ideal generated by $\leq 3$ quadratic forms
(the ring-theoretical classification was obtained in [E-K]). 
It is only in 5 variables that we can obtain non-finitely
generated Ext-algebras (only one case, just studied above in section 1) or Ext-algebras with
 an irrational Hilbert series (three cases).
 These three cases are as
follows (the numbering of cases is from [Ro4]) (in all these three cases the Hilbert series 
$R(t) = 1+5t+7t^2$):

\noindent Case 12:
$$
R_{12} = {E(x,y,z,u,v)\over (xy,xz+yu+zv,uv)} \quad{\rm with}
 \quad {1\over R_{12}^!(t)}=(1-2t)^2\prod_{n=1}^\infty (1-t^n)
$$
\noindent Case 20:
$$
R_{20} = {E(x,y,z,u,v)\over(yz+xu,yu+xv,zu+yv)}\quad{\rm with}
 \quad {1\over R_{20}^!(t)}=\prod_{n=1}^\infty (1-t^{2n-1})^5(1-t^{2n})^3
$$
\noindent Case 15:
$$
R_{15} = {E(x,y,z,u,v)\over(yz+xu,xv,zu+yv)}\quad{\rm with}
 \quad {1\over R_{15}^!(t)}=(1-2t)\prod_{n=1}^\infty(1-t^{2n-1})^3(1-t^{2n})^2
$$
But we can not see how any of these algebras could come up from some hyperplane arrangement.
If we study quotients of $E(x,y,z,u,v)$ with {\it four} quadratic forms there are still three
other quotients (this time with Hilbert series $R(t)=1+5t+6t^2$) which {\it might} have
irrational $R^!(t)$ :

\noindent Case 21:
$$\displaylines{
R_{21} = {E(x,y,z,u,v)\over (yz+xu,yu+xv,zu+yv,uv)} \quad{\rm with}
 \quad R_{21}^!(t)=1+5t+19t^2+65t^3+211t^4+667t^5+ \cr
+2081t^6+6449t^7+19919t^8+61425t^9+189273t^{10}+583008t^{11}+1795509t^{12}+\cr
+5529263t^{13}+17026752t^{14}+52431180t^{15}+161452384t^{16}+497162060t^{17}+\cr
+1530914456t^{18}+4714152439t^{19}+14516309322t^{20}+44700127353t^{21}+\cr
+137645268696t^{22}+423851580822t^{23}+ \ldots \cr
}$$
\noindent Case 22:
$$\displaylines{
R_{22} = {E(x,y,z,u,v)\over (yz+xu,yu+xv,zu+yv,zv)}\quad{\rm with}
 \quad R_{22}^!(t)=1+5t+19t^2+65t^3+211t^4+666t^5+ \cr
+2071t^6+6387t^7+19609t^8+60054t^9+183672t^{10}+561340t^{11}+ \cr
+1714894t^{12}+5237883t^{13}+15996477t^{14}+\ldots \cr
}$$
\noindent Case 33:
$$\displaylines{
R_{33} = {E(x,y,t,u,v)\over(yz+xu,xv,zu+yv,uv)}\,{\rm with}
 \, R_{33}^!(t)=1+5t+19t^2+65t^3+212t^4+675t^5+2125t^6+\cr
+6653t^7+\ldots (21\, {\rm terms}\,) +483131948638003t^{29}+1505474194810058t^{30}+ \ldots \cr
{\rm but\,\, the \,\, following \,\, formula \,\, gives \,\, in \,\, this \,\, last
\,\, case \,\, an \,\, indication
 \,\, about \,\, theta \,\, functions:}\cr
{1\over (1-t)^2 R^!_{33}(t)}= 1- 3t - t^2 + t^3 + 2t^4 + 3t^5 + t^6 +  \cr
+t^7 - t^8 - t^9 - 2t^{10}- t^{11} - 3t^{12} - t^{13} - t^{14} - t^{15} + t^{17} +\cr
+ t^{18} + 2t^{19}+ t^{20} + t^{21} + 3t^{22} + t^{23} + t^{25} + t^{26} - t^{29} - t^{30} \ldots \cr
}$$
leading to the predictions that the coefficient for $t^{31}$ should be $-2$ and that more 
precisely the series
should continue as $-2t^{31}-t^{32}-t^{33}-t^{34}-3t^{35}+ \ldots$.
But using the Backelin et al. programme BERGMAN [B] we have for the moment only been able to 
calculate the preceding series in degrees $\leq 30$ and no precise theory is in sight.
But I still do not know if these last three cases come from some hyperplane arrangements.
In higher embedding dimensions (6,7,...) there are of course more irrational series
 and as we have indicated two of them in embedding dimension 7 come from complex hyperplane
 arrangements ...

\noindent {\it Remark 4.4.} The case $R^!_{20}$  (which comes from J\"urgen Wisliceny  and whose series was
determined up to degree 67 by Czaba Schneider [Sch], Theorem 6.1 ) 
was completely determined in the super-Lie
algebra case in [L\"o-Ro1] (where we had a periodicity 4). Here its treatment is easier
 (periodicity 2).
Note that we have only described above what happens in characteristic 0. In case 20 we have
different $R^!_{20}(z)$ in all characteristics and the same remark seems to be
applicable to the cases 21,22,33.

\mysec {5. Irrational or non finitely presented cases for other arrangements ?}

In the sections above we have found two classes of unexpected complex hyperplane arrangements.
An interesting question is to determine how rare those hyperplane arrangements are.
The simplest arrangements are the so-called graphic arrangements: we have a simple graph $\Gamma$
given with $n$ vertices and $t$ edges. The corresponding hyperplane arrangement
 ${\cal A}_{\Gamma}$ in ${\bf C}^n$ is defined by:
$$
{\cal A}_{\Gamma} = \{x_i-x_j\} \, {\rm where \,} i < j, \,{\rm and \, where}\, (i,j)\, {\rm \,is\,
 an\, edge\, of} \, {\Gamma}
$$
Such a graph leads as in section 1 to an Orlik-Solomon algebra which can be written
in the form $R_{\Gamma} \otimes E(z)$ where
$R_{\Gamma}$ is a quotient of the exterior algebra in
$t-1$ variables by homogeneous forms, and $E(z)$ is the exterior algebra in one variable z.
It is therefore sufficient to analyze $R_{\Gamma}$.
Now Lima-Filho and Schenk have recently proved [Li-Sch] that the Hilbert series of all $R_{\Gamma}^!$
are rational of a special form, and therefore it follows that the Hilbert series of
the Ext-algebra of the Orlik-Solomon algebra of ${\cal A}_{\Gamma}$ is always rational,
at least for those cases where the Hilbert series of $R_{\Gamma}$  has the cube of its
maximal ideal equal to 0 (and maybe in all cases, cf. remarks below).
 Indeed formula (2.17) above can be applied
and gives an explicit rational formula. But non-finitely presented 
Ext-algebras can indeed occur for some graphs:
In the book by about graphs by Harary [H] there is at the
end an explit list of simple graphs with $\leq 6$ vertices. We have gone through that list 
completely
and we can report part of the results as Theorem 5.1 below 
(note that it is sufficient to analyze connected
graphs, since the Orlik-Solomon algebra decomposes as a tensor product of the algebras
corresponding to the connected components of the graph).
It is known that the number of simple connected graphs with $n$ vertices increases
rapidly with $n$, according to the following table:
\halign{#\hfill\quad&&\hfill#\quad\cr
$n$ = Number of vertices of the graph: & 2& 3& 4& 5& 6& 7& 8 \cr
Number of simple connected graphs with $n$ vertices: & 1& 2& 6& 21& 112& 853& 11117\cr} 
We also use the numbering of the simple connected
graphs from the home page of of Brendan McKay,

{\tt http://cs.anu.edu.au/\~{}bdm/data/graphs.html}

There you can download e.g. a file called graph4c.6g giving all the connected graphs
on 4 vertices in so-called .g6-format, i.e.
this file is in a strange compressed form. It looks like
\medskip
{\small CF

CU

CV

C{}]

C\^\

C\~{}
}

\noindent Next you download the showg programme from the page

{ \tt http://cs.anu.edu.au/\~{}bdm/data/formats.html}

\noindent and then e.g. the command 

\noindent {\tt showg -eo1 graph4c.6g connected-graphs-of-order4}

\noindent produces the file connected-graphs-of-order4 which looks as follows:

{\small Graph 1, order 4.

4 3

1 4 $\,\,$ 2 4 $\,\,$ 3 4

Graph 2, order 4.

4 3

1 3 $\,\,$ 1 4 $\,\,$ 2 4

Graph 3, order 4.

4 4

1 3 $\,\,$ 1 4 $\,\,$ 2 4 $\,\,$ 3 4

Graph 4, order 4.

4 4

1 3 $\,\,$ 1 4 $\,\,$ 2 3 $\,\,$ 2 4

Graph 5, order 4.

4 5

1 3 $\,\,$ 1 4 $\,\,$ 2 3 $\,\,$ 2 4 $\,\,$ 3 4

Graph 6, order 4.

4 6

1 2 $\,\,$ 1 3 $\,\,$ 1 4 $\,\,$ 2 3 $\,\,$ 2 4 $\,\,$ 3 4}

\noindent Thus we see that there are 6 connected graphs with 4 vertices. The first line after graph 4
only says that there are 4 vertices and 4 edges.
The second line

1 3 $\,\,$ 1 4 $\,\,$ 2 3 $\,\,$ 2 4  

\noindent says that these four edges are exactly those that connect vertices 1 and 3,
 vertices 1 and 4, vertices 2 and 3 and vertices 2 and 4,
 i.e. that the hyperplane arrangement is defined by
$$
(x_1-x_3)(x_1-x_4)(x_2-x_3)(x_2-x_4)
$$
This is the only graphical arrangement for a graph with 4 vertices (the graph of a square)
where the Orlik-Solomon algebra
is not a Koszul algebra. The Orlik-Solomon algebra satisfies howver the condition $L_3$ 
(cf. Remark 2.3 in section 2) since $m^4=0$ and the Ext-algebra is finitely presented, 
and $Ext^*_R(k,k)$ has a rational Hilbert series.

For graphs of order 5 and 6 we have the following

{\bf THEOREM 5.1.}
a) Among the 21  connected graphs with 5 vertices, 15 give rize to Orlik-Solomon algebras 
(OS-algebras) that are Koszul.
Among the 6 remaining non-Koszul algebras only {\it one} (corresponding to Graph 19 = the graph of
 a pyramid with a square basis)
 gives rise to a hyperplane  arrangement where the  Ext-algebra of the OS-algebra is not
finitely presented, graphs 5,7,15,17,19 give OS-algebras that satisfies $L_3$,
and graph 14 (the graph of a pentagon) has an OS-algebra that satisfies $L_4$,
and since all $R(z)$ and $R^!(z)$ are rational, the Hilbert series of the 21 Ext-algebras 
are rational.

 b) Among the 112 connected graphs with 6 vertices, 34 give rize to Orlik-Solomon algebras
 (OS-algebras) that are Koszul.
Among the 78 remaining non-Koszul algebras only 7 (corresponding to graphs 71,74,100,102,107,108,109)
 have non-finitely presented Ext-algebras
of their
OS-algebras and one (corresponding to graph 98) has a finitely presented Ext-algebra, which however
 has an infinite number
of relations between the relations.
The condition $L_5$ is satisfied in one case ( the graph of a 6-gon; more generally $L_{n-1}$ is satisfied
for the graph of an n-gon).
The condition $L_4$ is satisfied for the graphs 48,95,98 and
the condition $L_3$ is satisfied for the graphs 

\noindent 11,13,25,33,36,39,42,44,46,51,53,57,61,63,66,68,72,73,81,87,92,99,100,102,106,107,108,109.

\noindent There are 5 graphs for which no condition $L_n$ is verified: 38,71,74,96,97. But also for these graphs the Hilbert series of the
Ext-algebra can be analyzed and proved to be rational, so all these 112 graphs give rational series.

SKETCH OF PROOF OF PART OF THE THEOREM: a) The case of the graph 19 with 5 vertices gives rize to
 the arrangement defined by the polynomial
$$
(x_1-x_3)(x_1-x_4)(x_1-x_5)(x_2-x_3)(x_2-x_4)(x_2-x_5)(x_3-x_5)(x_4-x_5)
$$  
corresponding to a pyramid, where the vertex $x_5$ is at the top of the pyramid, and $x_1,x_3,x_2,x_4$ are at the basis.
We have 8 factors and the OS-algebra is in 8 variables:
$$
{E(e_1,e_2,e_3,e_4,e_5,e_6,e_7,e_8) \over {\scriptscriptstyle ((e_1-e_7)(e_3-e_7),(e_2-e_8)(e_3-e_8),(e_4-e_7)(e_6-e_7),(e_5-e_8)(e_6-e_8),
(e_1-e_2)(e_2-e_5)(e_4-e_5))}} \leqno(5.1)
$$
Let us now introduce new variables $x_i=e_i-e_5$ for $i\neq 5$ and $z=x_5$. Our algebra (5.1) becomes as earlier the tensor product of the
exterior algebra $E(z)$ and the algebra in seven variables ($x_5$ is missing !):
$$
R = {E(x_1,x_2,x_3,x_4,x_6,x_7,x_8) \over (x_1x_3-x_1x_7+x_3x_7,x_2x_3-x_2x_8+x_3x_8,x_4x_6-x_4x_7+x_6x_7,x_6x_8,x_1x_2x_4)}
$$
It is now easy to see that the annihilator of $x_6$ in $R$ is generated by $x_6$ and $x_8$ and and
 similarly that the annihilator of $x_8$
is generated by $x_6$ and $x_8$. Furthermore the intersection of the two ideals $x_6$ and $x_8$ is 0.
 Thus the ideal $a=(x_6,x_8)$ is a direct sum and $S=R/a$ has a linear resolution over $R$; 
more precisely we have $P_R^S(x,y) = 1/(1-2xy)$. 
 Now apply a result by Rikard B{\o}gvad [B{\o}]
 which says that if $R\rightarrow S=R/a$ is an algebra map such that
$R/a$ has a linear $R$-resolution, then the map $R \rightarrow S$ is a large map in the sense of 
Levin [Le3]. In the proof of Lemma 2.3 b), page 4 in [B{\o}], there is a slight
misprint on line 9 of the proof which should be:

" $\ldots \eta :R \rightarrow S$ is 1-linear, i.e. that $Tor_{i,j}^R(S,k)=0$ if $i\neq j \ldots$ "

\noindent This has the consequence that the
double series $P_R(x,y) = P_R^S(x,y)P_S(x,y)$ (i.e. the change of rings spectral sequence degenerates,
 cf. [Le3, Theorem 1.1, p. 209]. In our case
$P_R^S(x,y)=1/(1-2xy)$ and $S$ now becomes the quotient of an exterior algebra in five variables:
$$
S = {E(x_1,x_2,x_3,x_4,x_7) \over (x_1x_3-x_1x_7+x_3x_7,x_4x_7,x_2x_3,x_1x_2x_4)}
$$
But this algebra is essentially the algebra (1.2), treated in section 1 above. Indeed,
 let us first make the substitutions:
$ x_1 \longrightarrow x_1+x_3,\,\, x_7 \longrightarrow x_7+x_3$ and then interchange $x_2$ and $x_3$ and also interchange $x_1$ and $x_7$.
We get the same algebra as in  (1.2) in section 1; the only difference is that now the last variable
 is called $x_7$ (and not $x_5$ as in
section 1). It follows that the double series of $S$ is given by 
$$
{1 \over P_S(x,y)}= {1-6xy+12x^2y^2-x^2y^3-8x^3y^3\over 1-xy}
$$
so that
$$
{1 \over P_R(x,y)}= {(1-2xy)(1-6xy+12x^2y^2-x^2y^3-8x^3y^3)\over 1-xy } \leqno(5.2)
$$
From (5.2) we can now read off that $R(z)=1+7z+17z^2+14z^3=(1+2z)(1+5z+14z^3)$
 and that $R^!(z) = (1-z)/(1-2z)^4$ so that the formula $L_3$ holds true.
 One then proves that $gldim(R^!) = 4$ and that $Tor^{R^!}_{4,i}(k,k)$ is 1-dimensional for $i\geq 4$
and $Tor^{R^!}_{3,i}(k,k)=0$ for $i\neq 3$ so that 
 the algebra in section 1 which there needs an infinite number of generators 
now comes back here as a subtle part of a graphic arrangement
whose Ext-algebra is slightly better in that it is finitely generated but not finitely presented.
Note in particular that $Ext^*_R(k,k)$ is {\it not} the tensor product of $Ext^*_S(k,k)$
and the free algebra on two variables of degree 1. 
The other non-Koszul cases in Theorem 5.1 a) are simpler and treated in a similar way.
For graphs with 6 vertices (Theorem 5.1 b) similar procedures are used and the most complicated case is case 109, where
the arrangement is defined by the polynomial:
$$
\scriptstyle (x_1-x_3)(x_1-x_4)(x_1-x_5)(x_1-x_6)(x_2-x_3)(x_2-x_4)(x_2-x_5)(x_2-x_6)(x_3-x_5)(x_3-x_6)(x_4-x_5)(x_4-x_6)
$$ 
where we still have a non-finitely presented Ext-algebra.
But the most interesting pairs of examples are 107 and 74 which have the same Hilbert series {\it both} for the Orlik-Solomon algebra
{\it and} for the quadratic dual $R^!$ and the first algebra satisfies $L_3$ and the second does not. But they have both 
Ext-algebras that are not finitely presented. In particular the Tor (or Ext) of the Orlik-Solomon
 algebras differ, the first difference
occurs for $Tor^{\cal OS}_{4,6}(k,k)$ which has dimension 9 for the case 107 and 
dimension 10 for case 74.
A similar phenomenon occurs for the cases 87 and 71 where the first $Tor^{\cal OS}_{4,6}(k,k)$ has dimension 16 for the case 87
(but here the Ext-algebra is finitely presented) and the second case 71 has $Tor^{\cal OS}_{4,6}(k,k)$ of dimension 17 and the Ext-algebra
is not finitely presented in that case.

{\it Remark 5.2.} I have also studied many of the 853 cases corresponding to graphs with seven vertices.
Everything in sight leads to rational Hilbert series for the Ext-algebras.

{\it Remark 5.3.} When I lectured about this at Stockholm University, 
J\"orgen Backelin made an interesting observation:

1) For graphds with 5 vertices the only case of non-finitely presented Ext-algebras comes from the 
case when you remove two disjoints edges from the complete graph on 5 vertices (case 19).

2) If you remove three disjoint edges in the complete graph on 6 vertices you get the case 109 which
 is the most complicated one for graphs of order 6.

3) This leads to a conjecture that if you remove $[n/2]$ disjoint edges from
the complete graph on n vertices ($n\geq 7$) then you should get a very interesting situation.

{\it Remark 5.4.} Here is another example: 
The Example 1.3 of [Fi-Sch] where G is "the one skeleton of the Egyptian pyramid and the
 one skeleton of a tetrahedron sharing a single triangle" we have 
$1/R^!(t)=(1-2t)^4(1-3t)$ and $R(t) =1+11t+48t^2+103t^3+107t^4+42t^5$
leading to the formula
$$
{1\over P_R(x,y)} = 1-11\,x^2y^2+48\,x^3y^3-x^2y^3-104\,x^3y^3+5\,x^3y^4+112\,x^4y^4-6\,x^4y^5-
48\,x^5y^5
$$
i.e. this is another one of the cases where (2.17) is true but $m^3 \neq 0$.

Furthermore, the global dimension of $R^!$ is 5 and the Ext-algebra is finitely generated but
{\it not} finitely presented.

{\it Remark 5.5.} The preceding results (which should satisfy the request of one of the referees)
 show that the behaviour of the graphic arrangements in the
non-Koszul case are rather unpredictable (but the irrational case is rare and it is quite
probable that it does not occur for graphic arrangements). 

\mysec {6. Questions of Milnor, Grigorchuk, Zelmanov, de la Harpe and irrationality.} 

In section 4 we have presented two hyperplane arrangements, the Mac Lane arrangement $ML$ and
an easier variant $mleas$ =(4.1) which both have irrational Hilbert series for the 
corresponding $R^!$.
In section 4 we also presented the three possibilities for irrational Hilbert series for $R^!$
when $R$ is an arbitrary quotient of an exterior algebras in five variables with an
ideal generated by three quadratic forms (the only possibilities): Cases $R_{12},R_{20}$
 and $R_{15}$.
It is seems difficult to achieve similar examples for hyperplane arrangements using five variables.
 But $mleas$ can be considered as a higher variant of $R_{12}$
and similarly $ML$ can be considered as a higher variant of $R_{15}$ 
(indeed in this last case there are central
elements in odd degrees $\geq 3$ , so that we get the exponents 2,2,2,... in the infinite
product formula for case $R_{15}^!(t)$ in section 4 when we have divided out the center.
But so far I have not found any hyperplane arrangement corresponding (or similar) to the irrational
case $R_{20}$ in section 4. If such a hyperplane arrangement existed 
it would in particular lead to results about growth of groups and groups of finite width.
 Let me be more precise:
First recall that if ${\cal A}$ is {\it any} finite complex hyperplane arrangement in ${\bf C}^n$
and if $G=G({\cal A})$ is the fundamental group of the complement of the union of the corresponding
hyperplanes in ${\bf C}^n$ then $G$ is finitely presented --- indeed the complement
 has the homotopy type of a finite CW-complex ([Or], Proposition 5.1, p. 43 ). Let
$$
G=G_1 \supseteq G_2 \supseteq G_3 \supseteq \ldots 
$$ 
be the descending lower central series of $G$ defined inductively by $G_1 = G$ and
$G_k=[G_{k-1},G_1]$ (for $k \geq 2$). We have a structure of a graded Lie-ring 
(which can have torsion):
$$
gr(G) = \bigoplus_{i\geq 1} {G_i \over G_{i+1}}
$$
where the graded Lie structure is defined as follows: Let $\bar x$  and $\bar y$
be elements in $G_i/G_{i+1}$ and $G_j/G_{j+1}$ respectively,and let them be represented by $x$ and $y$
in $G_i$ and $G_j$. Then $xyx^{-1}y^{-1}$ lies in $G_{i+j}$ and its image in
$G_{i+j}/G_{i+j+1}$ is denoted by $[\bar x,\bar y]$.
It was proved by Kohno [K] that we have an isomomorphism of graded Lie algebras
$$
gr(G) \otimes_Z {\bf Q}   \simeq  \bigoplus \eta^i \leqno(6.1)
$$
where $\eta$ is the Lie algebra of primitive elements in the subalgebra generated by \break
$Ext^1_{OS({\cal A})}(Q,Q)$ of the Yoneda Ext-algebra of the Orlik-Solomon algebra of ${\cal A}$.
Now $G_{ML}$ and is $G_{mleas}$ are finitely presented groups.
Therefore, if we could find a hyperplane arrangement (probably in high embedding dimension) 
corresponding to the case $R_{20}$ or similar,
we would have at the same time found a {\it finitely presented group} $G$ such that the
Lie algebra $gr(G) \otimes_Z {\bf Q}$ is infinite and of finite width (i.e the dimensions of the
$\eta^i$ in (6.1) are bounded
(for further terminology and results I refer to
the surveys [delH],[Ba-Gr], [Gr-P] and the literature cited there). If so one would probably be
close to {\it finitely presented} groups having intermediate growth. Note that it is not expected that
such groups exist (cf.  Conjecture 11.3 in [Gr-P], where one states two lines earlier
that the existence of such groups is a major open problem in the field,  
 and research problem VI.63 on page 297 of [delH]).
But our two groups corresponding to the Mac Lane arrangement $ML$ and its easier variant $mleas$
give at least {\it finitely presented groups} with irrational growth series (Hilbert series)
of $U(gr(G)\otimes_Z {\bf Q})$.

\mysec {7. Final remarks.}
It is interesting to note that 32-33 years ago Jean-Michel Lemaire was in his thesis
[Lem1] inspired by the Stallings
group-theoretical example [St] (now used again in [Di-Pa-Su] !) to construct a finite simply-connected CW-complex $X$ such that
the homology algebra of the loop space $H_*(\Omega X,{\bf Q})$ was not finitely presented (not
even finitely generated). In [Ro1] I used a general recipe which in particular could be used to
translate Lemaire's results to local commutative ring theory 
to obtain a local ring $(R,m)$ such that the Yoneda  Ext-algebra ${\rm Ext}^*_R(k,k)$
was not finitely generated, thereby solving in the negative a problem by Gerson Levin [Le1].
The example in section 1 above is just a skew-commutative variant of my example in [Ro1],
 but with a quick direct proof, which hopefully should satisfy mathematicians working with
 arrangements of hyperplanes. The theory of section 2 above, combined with more difficult
variants of the later developments in the 1980:s about a question of Serre-Kaplansky [Lem2] 
are here shown to be useful for solving the second problem of Denham-Suciu [De-Su]. 

\bigskip
REFERENCES
\bigskip
[An] Anick, David J., A counterexample to a conjecture of Serre. Ann. of Math. (2) 115 (1982), no. 1, 
1--33. and Comment: "A counterexample to a conjecture of Serre". 
Ann. of Math. (2) 116 (1982), no. 3, 661.

[B] Backelin, J\"orgen et al, BERGMAN, a programme for non-commutative Gr\"obner basis
calculations available at 
{\tt http://servus.math.su.se/bergman/}

[Ba-Gr] Bartholdi, Laurent; Grigorchuk, Rotislav, Lie methods in growth of groups
 and groups of finite width, London Math. Soc. Lecture Notes, p. 1-27, vol. 275, Cambridge 
Univ. Press, Cambridge 2000.

[B{\o}] B{\o}gvad, Rikard, Some homogeneous coordinate rings that are Koszul algebras,

\noindent {\tt arXiv:math.AG/9501011} v2, 25 sept 1995.

[delH] de la Harpe, Pierre, Topics in Geometric Group Theory, The University of Chicago Press,
Chicago and London, 2000.

[De-Su] Denham, Graham and Suciu, Alexander I, On the homotopy Lie algebra of an arrangement,
Michigan Mathematical Journal 54, (2006), no 2, 319-340.

[Di-Pa-Su] Dimca, Alexandru; Papadima, Stefan and Suciu, Alexander I, Non-finitenes properties of
fundamental groups of smooth projective varieties. Preprint available at 

\noindent {\tt arXiv:math.AG/0609456}

[E-K] Eisenbud, David; Koh, Jee, Nets of alternating matrices and the linear syzygy conjectures.
 Adv. Math. 106 (1994), no. 1, 1--35.

[Gr-P] Grigorchuk, Rostislav; Pak, Igor, Groups of Intermediate Growth: an Introduction for Beginners.

  {\tt arXiv:math.GR/0607384, 17 july 2006}.

[H] Harary, Frank, Graph theory. Addison-Wesley Publishing Co., Reading,
 Mass.-Menlo Park, Calif.-London 1969 ix+274 pp.

[K]  Kohno, T.,On the holonomy Lie algebra and the nilpotent completion of the fundamental
group of the complement of hypersurfaces, Nagoya Math. J. 92 (1983), 21--37.

[Lem1] Lemaire, Jean-Michel, Alg\`ebres connexes et homologie des espaces de lacets.
 (French) Lecture Notes in Mathematics, Vol. 422. Springer-Verlag, Berlin-New York, 1974.
 xiv+134 pp. 

[Lem2] Lemaire, Jean-Michel Anneaux locaux et espaces de lacets \`a s\'eries de Poincar\'e
 irrationnelles (d'apr\`es Anick, Roos, etc...). 
(French) [Local rings and loop spaces with irrational Poincare series (after Anick, Roos, etc.)]
 Bourbaki Seminar, Vol. 1980/81, pp. 149--156, Lecture Notes in Math., 901,
 Springer, Berlin-New York, 1981.

[Le1] Levin, Gerson, Two conjectures in the homology of local rings. J. Algebra 30 (1974), 56--74.

[Le2] Levin, Gerson, Modules and Golod homomorphisms. 
J. Pure Appl. Algebra 38 (1985), no. 2-3, 299--304.

[Le3] Levin, Gerson, Large homomorphisms of local rings, Math. Scand. 46 (1980), 209--215.

[Li-Sch] Lima-Filho, Paulo ; Schenck, Hal, Holonomy Lie algebras and the LCS-formula for
subarrangements of $A_n$, (2006), preprint available at 

{\tt http://www.math.tamu.edu/\~{}schenck/glcs.pdf}

[L\"o1] L\"ofwall, Clas On the subalgebra generated by the one-dimensional elements in the Yoneda 
Ext-algebra. Algebra, algebraic topology and their interactions (Stockholm, 1983), 291--338, 
Lecture Notes in Math., 1183, Springer, Berlin, 1986.

[L\"o2] L\"ofwall, Clas, On the homotopy Lie algebra of a local ring. J. Pure Appl. Algebra 38 (1985),
 no. 2-3, 305--312.

[L\"o3] L\"ofwall, Clas, {\tt liedim.m} , a Mathematica programme to calculate (among other things)
the ranks of a finitely
presented graded Lie algebra (2001), available at

{\tt http://www.math.su.se/\~{}clas/liedim/ }

[L\"o4] L\"ofwall, Clas, Appendix B to [Ro3]

[L\"o-Ro1] L\"ofwall, Clas; Roos, Jan-Erik, A nonnilpotent $1$-$2$-presented graded Hopf algebra whose Hilbert
 series
converges in the unit circle. Adv. Math. 130 (1997), no. 2, 161--200. 

[L\"o-Ro2] L\"ofwall, Clas; Roos, Jan-Erik Cohomologie des alg\`ebres de Lie 
gradu\'ees et s\'eries de Poincar\'e-Betti non rationnelles. (French)
 C. R. Acad. Sci. Paris S\'er. A-B 290 (1980), no. 16, A733--A736.

[Mi-Mo] Milnor, John W.; Moore, John C., On the structure of Hopf algebras. 
Ann. of Math. (2) 81 (1965) 211--264. 

 [Or] Orlik, Peter, Introduction to arrangements. CBMS Regional Conference Series in Mathematics, 72. 
Published for the Conference Board of the Mathematical Sciences, Washington, DC;
 by the American Mathematical Society, Providence, RI, 1989. x+110 pp.

[Pa-Su] Papadima, Stefan and Suciu, Alexander I., When does the associated graded Lie algebra of an
arrangement group decompose ?, Comment. Math. Helv. 81 (2006), 859-875.

[Ro1] Roos, Jan-Erik, Relations between Poincar\'e-Betti series of loop spaces
 and of local rings. S\'eminaire d'Alg\`ebre Paul Dubreil 31\`eme ann\'ee (Paris, 1977--1978),
 pp. 285--322, Lecture Notes in Math., 740, Springer, Berlin, 1979.

[Ro2] Roos, Jan-Erik, On the use of graded Lie algebras in the theory of local rings.
 Commutative algebra: Durham 1981 (Durham, 1981), pp. 204--230, 
London Math. Soc. Lecture Note Ser., 72, Cambridge Univ. Press, Cambridge-New York, 1982. 

[Ro3] Roos, Jan-Erik, A computer-aided study of the graded Lie algebra of a local commutative Noetherian
 ring. J. Pure Appl. Algebra 91 (1994), no. 1-3, 255--315.

[Ro4] Roos, Jan-Erik, Homological properties of quotients of exterior algebras,
in preparation; cf. Abstracts Amer.Math.Soc., 21 (2000), 50-51.

[Ro5] Roos, Jan-Erik, On computer-assisted research in homological algebra, Mathematics and Computers
 in Simulation, 42, (1996), 475-490.

[Ro6] Roos, Jan-Erik, Paper about the irrationality of hyperplane arrangements and
growth of groups, in preparation.

[Ro6] Roos, Jan-Erik, Homology of loop spaces and of local rings.
 18th Scandinavian Congress of Mathematicians (Aarhus, 1980), pp. 441--468, 
 Progr. Math., 11, Birkh\"auser, Boston, Mass., 1981.

[Sch] Schneider, Czaba ,Computing Nilpotent Quotients in Finitely Presented Lie \break Rings.
 Discrete Mathematics \& Theoretical Computer Science. 1(1), pages 1-16, 1997.

[Sh-Yu] Shelton, Brad; Yuzvinsky, Sergey, Koszul algebras from graphs
 and hyperplane arrangements. J. London Math. Soc. (2) 56 (1997), no. 3, 477--490.

[St] Stallings, John ,A finitely presented group whose 3-dimensional integral homology 
is not finitely generated. Amer. J. Math. 85 (1963) 541--543.

[Su] Suciu, Alexander I. Fundamental groups of line arrangements: enumerative aspects.
 Advances in algebraic geometry motivated by physics (Lowell, MA, 2000), 43--79,
 Contemp. Math., 276, Amer. Math. Soc., Providence, RI, 2001.

\par}

\end